\documentclass[12pt]{amsart}
\usepackage{amssymb}
\usepackage[colorlinks=true,linkcolor=black,citecolor=black,urlcolor=black]{hyperref}
\hypersetup{pdftitle={Stability conditions for families of holomorphic endomorphisms on projective spaces}}
\usepackage[a4paper,width=16.2cm,top=3cm,bottom=3cm]{geometry}

\theoremstyle{plain}
\newtheorem{thm}{Theorem}[section]
\newtheorem{prop}[thm]{Proposition}
\newtheorem{cor}[thm]{Corollary}
\theoremstyle{definition}
\newtheorem{rem}[thm]{Remark}

\theoremstyle{plain}

\numberwithin{equation}{section}
\newcommand{\C}{\mathbb C}
\renewcommand{\P}{\mathbb P}
\newcommand{\F}{\mathcal F}
\newcommand{\J}{\mathcal J}
\newcommand{\Lcal}{\mathcal L}
\newcommand{\supp}{\operatorname{supp}}

\usepackage{color}
\usepackage{xcolor}

\definecolor{revisionpurple}{RGB}{128,0,128}

\begin{document}
\pdfoutput=1
\title[Equilibrium webs and stability for families of endomorphisms]{Equilibrium webs and 
stability 
 for 
arbitrary families of endomorphisms of projective spaces}

\author[F.~Bianchi]{Fabrizio Bianchi}
\address{Dipartimento di Matematica, Università di Pisa, Largo Bruno Pontecorvo 5, 56127 Pisa, Italy}
\email{fabrizio.bianchi@unipi.it}

	\author{Zhuchao Ji}

\address{Institute for Theoretical Sciences, Westlake University, Hangzhou 310030, China}

\email{jizhuchao@westlake.edu.cn}

\author{Junyi Xie}

\address{Beijing International Center for Mathematical Research, Peking University, Beijing 100871, China}

\email{xiejunyi@bicmr.pku.edu.cn}
\date{\today}

\begin{abstract}
We consider arbitrary 
holomorphic families of endomorphisms of
a complex projective space.
We prove that the existence of an equilibrium web implies the
holomorphic motion of all repelling 
cycles contained in the Julia set, i.e., the support of the unique measure of maximal entropy.
 Together with previous 
results of Berteloot, Bianchi,
Br\'evard, Buff, and Dupont, this
completes the extension of the stability-bifurcation theory to
 arbitrary
families of holomorphic endomorphisms.
The same conclusions also 
hold for holomorphic families of polynomial-like
maps of large topological degree.

\end{abstract}
\maketitle

\section{Introduction}

\subsection{The one-dimensional theory}

The stability and bifurcation theory of rational maps was developed in the early 1980s in the works of   
Lyubich \cite{Lyu83} and Mañé, Sad and Sullivan \cite{MSS83}, see also 
\cite{DH84,Lev81,McM94,McM00}. 
These works established several equivalent characterizations of stability. In particular, stability can be expressed in terms of the holomorphic motion of the Julia sets,
 the holomorphic motion of the repelling cycles, the normality of the critical orbits, or the absence of non-persistent neutral cycles. They also proved the fundamental fact that the stability locus is open and dense in every holomorphic family. 

The classical arguments rely on several features which are specific to one complex dimension. Repelling periodic points are dense in the Julia set and move locally holomorphically with the parameter. The $\lambda$-lemma allows one to extend their motion to the full Julia set, while Hurwitz's theorem prevents distinct moving points from colliding. Moreover, the critical set is finite, and the Picard--Montel theorem relates stability to the normality of finitely many critical orbits.

A different approach to bifurcations, based on potential theory, was introduced by DeMarco \cite{DeM01,DeM03}. 
If $\mu_\lambda$ denotes the equilibrium measure of $f_\lambda$
\cite{FLM83,Lyu83e}
 and
\[
 L(\lambda):=\int_{\mathbb P^1}\log|f_\lambda'|\,d\mu_\lambda
\]
is its Lyapunov exponent, then
the function $\lambda \mapsto L(\lambda)$
 is plurisubharmonic and the support of the positive closed current
$T_{\mathrm{bif}}:=dd^cL$
coincides with the bifurcation locus.
Beyond detecting the bifurcation locus, the current
$T_{\mathrm{bif}}$ and its powers arise as limits of suitably
normalized currents of integration over dynamically defined parameter
loci, such as critical-orbit relations and prescribed-multiplier
conditions; see \cite{DF08,BB11, BE09, GOV19} and the references therein.
We refer to the surveys
\cite{Ber13,Duj14} for an
introduction
 to the one-dimensional bifurcation theory and its
pluripotential approach.
These currents and measures also provide an important bridge
with arithmetic dynamics, through arithmetic equidistribution and
height theory, with applications to unlikely-intersection problems
and to the geometry of postcritically finite loci;
 see, for instance,
\cite{BD13,FG22,JX23,DMY26,DeM18} and the references therein.

\subsection{Higher-dimensional stability and equilibrium webs}

For holomorphic families of endomorphisms of 
 $\mathbb P^k = \mathbb P^k (\mathbb C)$,
 the classical one-dimensional proof cannot be directly adapted. The critical set has positive dimension, there is no analogue of the Picard--Montel argument involving finitely many critical orbits, and no general higher-dimensional $\lambda$-lemma is available. More fundamentally, Hurwitz's theorem no longer prevents two distinct holomorphic graphs with values in $\mathbb P^k$ from meeting. Thus, even when many points in the Julia sets can be followed holomorphically, the resulting motion may be branched.

To overcome these difficulties, Berteloot, 
Dupont,
and the first named author
 introduced in \cite{BBD18} the notions of 
\emph{equilibrium web} and 
\emph{equilibrium lamination}. 
Their approach also builds on the
potential-theoretic theory of bifurcation currents, previously
extended to higher-dimensional projective endomorphisms by
Bassanelli and Berteloot \cite{BB07} and Pham \cite{Pha05}.
For general background on
higher-dimensional holomorphic dynamics, we refer to
\cite{FS94,DS10}.
We 
also refer to \cite{AB23} for an equidistribution theorem in
parameter space and to \cite{GTV26}, where bifurcation measures 
\cite{AB22}
play
a key role in proving that postcritically finite endomorphisms of
$\mathbb P^k$ are not Zariski dense when $k\geq2$; see also \cite{GT26}.
Recently, building on the higher-dimensional bifurcation theory \cite{BBD18, BB26}, 
the second and third named authors \cite{JX26} establish genus and
gonality growth for dynatomic curves in one-parameter
families of endomorphisms of $\mathbb P^k$, with applications to
the
 dynamical uniform boundedness conjecture.

An equilibrium web can be viewed as a possibly branched measurable holomorphic motion of the equilibrium measures, while an equilibrium lamination gives an unbranched motion on a set of full equilibrium measure. Using ergodic and pluripotential methods, they obtained a higher-dimensional counterpart of the Lyubich--Mañé--Sad--Sullivan theory. In particular, they related the holomorphic motion of repelling $J$-cycles
(i.e., the repelling cycles contained in the Julia set), 
the pluriharmonicity of the sum of the Lyapunov exponents, the existence of equilibrium webs and laminations, and the stability of the postcritical dynamics.
We refer to the survey \cite{BB18} for an informal presentation of this theory,
and for a presentation of the one-dimensional theory with this point of view.
An important difference with the one-dimensional picture is that
stability is no longer dense in general: in higher dimensions the
bifurcation locus may have non-empty interior. We refer to
\cite{BT17,Duj17,Taf21,Bie19} for examples and different mechanisms
producing such robust bifurcations.
We also refer to
\cite{DL15,BD17}
for
a parallel theory of stability
and bifurcations for polynomial automorphisms of $\mathbb C^2$.

A key distinction in \cite{BBD18} is that between general and
\emph{acritical} equilibrium webs. The latter give no mass to graphs
meeting the grand orbit of the critical set. For arbitrary families,
the holomorphic motion of the repelling $J$-cycles gives an
acritical equilibrium web \cite[Corollary~1.7]{BBD18}. Conversely,
after replacing an acritical equilibrium web by an ergodic one using
\cite[Proposition~2.4]{BBD18}, the local construction in the proof of
\cite[Theorem~4.1]{BBD18} produces an equilibrium lamination.
Together with \cite[Proposition~1.2]{BB26}, these results show that
the holomorphic motion of the repelling $J$-cycles, the existence of
an acritical equilibrium web, and the existence of an equilibrium
lamination are equivalent for arbitrary connected families.

What remained open for arbitrary families was whether the
acriticality assumption could be removed. In
\cite[Proposition~5.6]{BBD18}, the mere existence of an equilibrium
web was shown to imply the holomorphic motion of all repelling
$J$-cycles only when $k=2$ or for open subsets of the full parameter
space. Consequently, the full equivalence in \cite{BBD18} was
established for simply connected open subsets of the full parameter
space of endomorphisms of $\mathbb P^k$, and for arbitrary simply
connected families when $k=2$. Most of the other arguments worked
for arbitrary families in every dimension.

 We observe here that the
 simple connectedness assumption has a different, purely global
origin: it rules out monodromy. This obstruction already occurs in
dimension one. Consider the quadratic family $z\mapsto z^2+c$ over
the complement of the Mandelbrot set. The family is locally stable,
but the two fixed points
$(1\pm\sqrt{1-4c})/2$
are exchanged by analytic continuation along a loop surrounding the
Mandelbrot set. Thus the repelling cycles do not admit a global
holomorphic motion over this parameter space.

To deduce the holomorphic motion of all repelling $J$-cycles from the
existence of an equilibrium web in the cases considered in
\cite{BBD18}, the authors studied what happens when a repelling
periodic point ceases to be repelling.
The goal was to produce a 
\emph{virtually repelling}
 Siegel periodic point and then contradict the existence of an equilibrium web. This required a detailed analysis of the multipliers, Brjuno linearization, and the possible persistent resonances and non-diagonalizability. Such exceptional configurations can be excluded by perturbing inside the full parameter space. They can also be avoided in dimension $2$, where, 
 at the 
  boundary parameters
under consideration, 
  there is only one neutral direction. 
 For an arbitrary family in dimension at least $3$, however, this method did not give the desired implication.

Several subsequent works clarified and strengthened different parts of this picture. 
The first named author
 extended some of the stability theory to polynomial-like maps of large topological degree and introduced a weaker version of the motion of the repelling points
\cite{Bia19},
which was
 proved to be equivalent to the other conditions for any family.
Berteloot proved an all-or-nothing phenomenon for the bifurcation of repelling cycles \cite{Ber18}. 
Berteloot and Brévard recast postcritical normality in terms of a
ramification current and proved that its convergence domain gives a
local characterization of stability \cite{BBr23}.
Rakhimov and the first named author showed that, in a stable family,
the measurable holomorphic motion is unbranched almost everywhere
also with respect to much more general ergodic invariant measures
supported on the Julia set, with positive Lyapunov exponents and
giving no mass to the postcritical set \cite{BR24}.
  Further developments for stable families include the
equidistribution of weighted periodic graphs \cite{BiBr26} and the
propagation of equilibrium states \cite{BrR24}.
More recently, Berteloot and Buff proved that, for arbitrary simply
connected families, the holomorphic motion of a positive proportion
of the repelling $J$-cycles already implies the holomorphic motion
of all repelling $J$-cycles \cite{BB26}.
  Their result also incorporates the chain-reaction phenomenon of \cite{Ber18}.

The main point of the present paper is that no acriticality assumption
 on an equilibrium web
is needed 
 to deduce the other stability conditions. 
More precisely, for an arbitrary holomorphic family of
endomorphisms of $\mathbb P^k$, we prove the
 missing
 implication
\[
 \text{existence of an equilibrium web}
 \,\,\Longrightarrow\,\,
 \text{holomorphic motion of all repelling $J$-cycles}.
\]
Combined with the results mentioned above, this completes the extension to arbitrary families of the 
higher-dimensional counterpart of the classical stability theory of Levin, Lyubich, Mañé--Sad--Sullivan and DeMarco.

\subsection{Main results}

Our main result is the following.

\begin{thm}\label{thm:main}
Let $M$ be a connected complex manifold and let
$ f:M\times\mathbb P^k\longrightarrow M\times\mathbb P^k$
be a holomorphic family of endomorphisms of algebraic degree
$d\geq2$. If $f$ admits an equilibrium web over $M$, then the
repelling $J$-cycles move holomorphically over $M$.
\end{thm}

 Combining this theorem with the results of
\cite{BBD18,BBr23,BB26} gives the following higher-dimensional
stability theorem for arbitrary connected families; simple
connectedness is needed only for the final equivalence below.
Let $C_f$ be the critical divisor of $f$, counted with multiplicity,
and let $[C_f]$ be its current of integration.
For a positive
$(1,1)$-current $T$ on $M\times\P^k$ and a relatively compact open
set $W\Subset M\times\P^k$, denote by $\|T\|_W$ the mass of $T$ on
$W$. We denote by
$\pi_M:M\times\P^k\to M$
the canonical projection.
The precise definitions of the various conditions are given in Section \ref{s:defs}.

\begin{thm}\label{thm}
Let $M$ be a connected complex manifold and let
$f:M\times\P^k\to M\times\P^k$ be a holomorphic family of
endomorphisms of algebraic degree $d\geq2$. Consider the following
conditions:

\begin{enumerate}
\item[\textup{(A)}]
The repelling $J$-cycles move holomorphically over $M$.

\item[\textup{(A')}] 
A positive proportion of the
repelling $J$-cycles
 move holomorphically over $M$. 

\item[\textup{(B)}]
The function $L$ is pluriharmonic on $M$.

\item[\textup{(C)}]
The family $f$ admits an equilibrium web.

\item[\textup{(C')}]
The family $f$ admits an acritical equilibrium web.

\item[\textup{(D)}]
The family $f$ admits an equilibrium lamination.

\item[\textup{(E)}]
Every $\lambda_0\in M$ has a relatively compact neighborhood
$U$ such that
\[
\liminf_{n\to\infty}d^{-kn}
\bigl\|(f^n)_*[C_f]\bigr\|_{U\times\mathbb P^k}=0.
\]

\item[\textup{(E')}]
For every relatively compact open set $U\Subset M$,
\[
\bigl\|(f^n)_*[C_f]\bigr\|_{U\times\mathbb P^k}
=O\bigl(d^{(k-1)n}\bigr).
\]

\item[\textup{(F)}]
The convergence domain
$\Omega(R_f)$
 of the ramification current
$R_f$
 projects onto 
 the
whole parameter space $M$, i.e., we have $\pi_M\bigl(\Omega(R_f)\bigr)=M$.

\item[\textup{(F')}]
The ramification current converges locally on the whole
parameter--phase space, i.e., we have 
$\Omega(R_f)=M\times\P^k$.

\item[\textup{(G)}]
There are no Misiurewicz parameters in $M$.

\end{enumerate}
Then the following assertions hold:

\begin{enumerate}
\item[\textup{(i)}]
Conditions \textup{(A)}, \textup{(A')}, \textup{(C)},
\textup{(C')} and \textup{(D)} are equivalent.

\item[\textup{(ii)}]
Conditions \textup{(B)}, \textup{(E)}, \textup{(E')},
 \textup{(F)}, \textup{(F')} and \textup{(G)}
 are equivalent.

\item[\textup{(iii)}]
Each of the equivalent conditions in \textup{(i)} implies each of
the equivalent conditions in \textup{(ii)}.

\item[\textup{(iv)}]
If, in addition, $M$ is simply connected, then all the conditions
above are equivalent.

\item[\textup{(v)}]
If the equivalent conditions in \textup{(i)} hold, the equilibrium
web $\nu^{\rm eq}$ is unique. Moreover,
\[
\nu^{\rm eq}(\Lcal_1\mathbin{\triangle}\Lcal_2)=0
\]
for any two equilibrium laminations $\Lcal_1$ and $\Lcal_2$ of $f$.
\end{enumerate}
\end{thm}

Theorems~\ref{thm:main} and~\ref{thm} remain valid for holomorphic
families of polynomial-like maps of large topological degree
\cite{DS03}; see Remark~\ref{rem:polynomial-like}.
We also refer to \cite{BR24} for a local analogue involving webs and
laminations associated with more general ergodic invariant measures
supported on the Julia set, with positive Lyapunov exponents and
giving no mass to the postcritical set. Under the hypotheses specified
there, the existence of these objects is equivalent to the stability
conditions above.

\subsection{Outline of the proof of Theorem \ref{thm:main}}

The key new ingredient is a simple uniqueness criterion for
holomorphic graphs, see Proposition~\ref{prop:uniqueness},
which
strengthens \cite[Lemma~2.5(3)]{BBD18}.
Let $K$ be a compact forward invariant family of holomorphic graphs
and let $x$ belong to the fiber over a parameter $\lambda_0$. We
prove that, if
\[
 \sup_{n\geq0}\bigl\|D(f_{\lambda_0}^n)_xv\bigr\|=+\infty
 \qquad\text{for every }v\in T_x\mathbb P^k\setminus\{0\},
\]
then at most one graph in $K$ passes through $x$. Indeed, if two
distinct graphs meet at $x$, the first non-zero coefficient of their
difference in the parameter direction is transformed under iteration
by $D(f_{\lambda_0}^n)_x$. The hypothesis forces this coefficient to
become unbounded along a subsequence, whereas the compactness of $K$
and Cauchy's estimates give a uniform bound.

Since all Lyapunov exponents of the equilibrium measure are positive
\cite{BD99}, this criterion applies at almost every point in any
fixed fiber of an equilibrium web. It follows that the web is
unbranched on a completely invariant set of full equilibrium measure,
and that inverse orbits in this set admit unique lifts to the support
of the web. We then combine this uniqueness property with the local
unstable-manifold argument of Berteloot and Buff \cite{BB26}. If a
moving repelling cycle ceased to be repelling, the lifted inverse
orbits of a positive-measure family of graphs would force their
values at a suitable parameter to belong to the forward orbit of a
proper unstable manifold. This is impossible because the equilibrium
measure
 is the
Monge--Ampère measure of a current with continuous potentials and
therefore
 gives no mass to pluripolar sets
\cite{DS10}.

Thus the proof avoids the analysis of neutral multipliers, resonances
and linearization used in \cite{BBD18}. In this sense, the uniqueness
criterion plays the role of an almost-everywhere higher-dimensional
substitute for the non-intersection consequence of Hurwitz's theorem.
It is also complementary to \cite{BR24}: there, stability is used to
construct unbranched motions for measures with positive Lyapunov
exponents, whereas here the almost-everywhere unbranching of a given
equilibrium web is used to deduce stability.

\medskip

\noindent\textbf{Acknowledgements.}  
The authors would like to thank Fran\c{c}ois Berteloot  for 
 useful  comments on an earlier
 version of the paper.
This 
work
 has received funding from
 the 
ANR project 
TIGerS/ANR-24-CE40-3604
and from the MIUR Excellence Department Project awarded to the Department of Mathematics, University of Pisa, CUP I57G22000700001.
The first named
author is affiliated to the GNSAGA group of INdAM. The second named author  is supported by the National Key R\&D Program of China (No.
2025YFA1018300), NSFC Grant (No.12401106), and ZPNSF grant (No.XHD24A0201).
The third named author  is supported by NSFC Grant (No.12271007) and by the Xplorer Prize from
the New Cornerstone Science Foundation.

\medskip

\noindent\textbf{AI disclosure.}
The present project originated in discussions, motivated by
a preliminary version of
\cite{JX26}, about whether the equilibrium web arising in their
argument 
 needed to
be acritical in order to deduce stability. The key new ingredient of this paper is Proposition \ref{prop:uniqueness}, whose main idea was suggested by \mbox{GPT-6} Astra Ultra. The authors wrote the proof and independently verified the statement and every step of the argument. All other statements and proofs in the paper were developed by the authors.

\section{Equilibrium webs and stability notions}
\label{s:defs}
Let $M$ be a connected complex manifold and let
\[
 f:M\times\P^k\longrightarrow M\times\P^k,
 \qquad (\lambda,z)\longmapsto(\lambda,f_\lambda(z)),
\]
be a holomorphic family of endomorphisms of algebraic degree $d\geq2$.
For every $\lambda\in M$, let $\mu_\lambda$ be the equilibrium measure
of $f_\lambda$ and let $J_\lambda:=\supp\mu_\lambda$ be its Julia set.
We refer to \cite[Section~1]{DS10} for the basic properties of these
measures. Let
\[
 L(\lambda):=\int_{\P^k}\log|\operatorname{Jac}f_\lambda|\,d\mu_\lambda
\]
be the sum of their Lyapunov exponents.
The function $L$ is plurisubharmonic by \cite{Pha05};
see also \cite{DS10} and \cite[Section~3.1]{BBD18}.

\medskip

Following \cite[Definition~1.3]{BBr23}, we consider the
parameter--phase-space version of the ramification current of
Dinh and Sibony, originally introduced in the setting of
polynomial-like maps; see \cite{DS03} and
\cite[Chapter~2]{DS10}. Namely, set
\[
R_{n,f}:=d^{-kn}(f^n)_*[f(C_f)],\qquad n\geq0,
\]
where $[f(C_f)]$ denotes the current of integration on the analytic
hypersurface $f(C_f)$, and formally set
\[
R_f:=\sum_{n\geq0}R_{n,f}.
\]
The convergence domain $\Omega(R_f)$ is the set of points
$x\in M\times\P^k$ admitting a relatively compact neighborhood $W$
such that
\[
\sum_{n\geq0}\|R_{n,f}\|_W<+\infty.
\]

\medskip

Let  $\mathcal O(M,\P^k)$ be the space of all holomorphic maps from $M$ to $\P^k$ with the topology of locally uniform convergence.  Denote
\[
 \J:=\{\gamma\in\mathcal O(M,\P^k):
       \gamma(\lambda)\in J_\lambda\text{ for every }\lambda\in M\}.
\]
For $\gamma\in\mathcal O(M,\P^k)$, let $\Gamma_\gamma$ be its graph.
Define the evaluation maps $e_\lambda$ and the induced dynamical system $\F$ by
\[
 e_\lambda(\gamma):=\gamma(\lambda),\qquad
 (\F\gamma)(\lambda):=f_\lambda(\gamma(\lambda)).
\]
The map $\F$ acts continuously on $\mathcal O(M,\P^k)$ and preserves
$\J$. An \emph{equilibrium web} is by definition a compactly supported Borel
probability measure $\nu$ on $\J$ such that for every $\lambda\in M$, 
we have
\[
 \F_*\nu=\nu\qquad
 \text{and}
 \qquad
 (e_\lambda)_*\nu=\mu_\lambda\quad.
\]
This notion was introduced in \cite[Definition~1.3]{BBD18}.
Let
\[
\operatorname{GO}(C_f)
:=\bigcup_{a,b\geq0}f^{-a}\bigl(f^b(C_f)\bigr)
\]
be the grand orbit of $C_f$.
An equilibrium web $\nu$ is called \emph{acritical} if
\[
\nu\bigl(\{\gamma\in\J:
\Gamma_\gamma\cap\operatorname{GO}(C_f)\neq\varnothing\}\bigr)=0;
\]
see \cite[Definition~2.3]{BBD18}.

\medskip

For $n\geq1$, let $R_n(\lambda)$ be the set of repelling fixed points
of $f_\lambda^n$ which belong to $J_\lambda$. 
The union $\bigcup_{n\geq1}R_n(\lambda)$ is dense in
$J_\lambda$ \cite{BD99}.
We say that the  
\emph{repelling $J$-cycles move holomorphically over $M$} if, for every
$n$, there is a finite set $\mathcal R_n\subset\J$ such that the  following 
evaluation maps are bijective
\[
 e_\lambda:\mathcal R_n\longrightarrow R_n(\lambda)
 \quad\text{for every }\lambda\in M.
\]
The restriction to repelling cycles contained in $J_\lambda$ is
essential. Unlike in dimension one, in higher dimension a repelling
periodic point need not belong to $J_\lambda$; examples are given in
\cite[Section~6, Example~2]{HP94} and
\cite[Section~4.1]{FS01}. See also
\cite[Section~2.1.3]{Bia16} for a family in which a fixed point
remains repelling but belongs to $J_\lambda$ at one parameter and
lies outside $J_\lambda$ at another.

 Following \cite{Ber18},
 we
  say that \emph{a positive proportion of the
   repelling $J$-cycles move holomorphically over $M$} if,  for every $n\geq1$, there is a finite
$\mathcal F$-invariant set $\mathcal P_n\subset\mathcal J$ such
that $e_\lambda$ maps $\mathcal P_n$ injectively into
$R_n(\lambda)$ for every $\lambda\in M$, and
\[
\limsup_{n\to\infty}d^{-kn}\#\mathcal P_n>0.
\]
 A stronger notion, requiring 
\emph{asymptotically all} repelling $J$-cycles  to move holomorphically
was introduced and used in \cite{Bia19}.

\medskip

 An \emph{equilibrium lamination} is a
relatively compact subset $\Lcal\subset\J$ such that distinct
elements have disjoint graphs,
$\mu_\lambda(e_\lambda(\Lcal))=1$ for every $\lambda\in M$, each
graph avoids $\operatorname{GO}(C_f)$, and
$\F:\Lcal\to\Lcal$ is $d^k$-to-one
\cite[Definition~1.4]{BBD18}.

\medskip

A parameter $\lambda_0$ is called \emph{Misiurewicz} if there are a
neighborhood $V$ of $\lambda_0$, a holomorphic section
$p:V\to\P^k$ of repelling $J$-periodic points of a fixed period,
an integer $n\geq1$, and an irreducible component $H$ of
$f^n(C_f)\cap(V\times\P^k)$ such that
\[
(\lambda_0,p(\lambda_0))\in H,
\qquad \Gamma_p\not\subset H.
\]
We use here the formulation of
\cite[Definition~1.1]{Bia19}; see also
\cite[Definition~1.5]{BBD18}.

\section{Uniqueness of graphs via unbounded derivatives}\label{sec:uniqueness}

Throughout this section, $M$ is a connected complex manifold.

\begin{prop}\label{prop:uniqueness}
Let $K\subset\mathcal O(M,\P^k)$ be compact and satisfy
$\F(K)\subset K$. Fix $\lambda_0\in M$ and $x\in\P^k$. Assume that
\begin{equation}\label{eq:expansion}
 \sup_{j\geq0}\bigl\|D(f_{\lambda_0}^j)_xv\bigr\|=+\infty
 \quad\text{for every }0\neq v\in T_x\P^k.
\end{equation}
Then there is at most one $\gamma\in K$ such that
$\gamma(\lambda_0)=x$.
\end{prop}

When $K$ is the support of an equilibrium web and $x$ is a repelling
periodic point, the uniqueness assertion in
Proposition~\ref{prop:uniqueness} already follows from the argument
in the proof of \cite[Lemma~2.5(3)]{BBD18}. The extension beyond
periodic points will be essential in Corollary~\ref{cor:lifting};
see Remark~\ref{rem:compare-BBD}.

\begin{proof}
Suppose that $\gamma,\eta\in K$ are distinct and satisfy
$\gamma(\lambda_0)=\eta(\lambda_0)=x$.
Restricting to a small holomorphic disk through $\lambda_0$ on
which their germs are different, we may assume that $M=\mathbb D$ and
$\lambda_0=0$.

Choose affine coordinates $\psi_x$ with $\psi_x(x)=0$, and set
\[
 a(t):=\psi_x(\gamma(t)),\quad b(t):=\psi_x(\eta(t)),
 \quad
\text{and}
 \quad \delta(t):=a(t)-b(t).
\]
Write
\[
 \delta(t)=t^r u+O(t^{r+1}),\qquad u\ne0,
\]
where $r\geq1$ is the first order at which the Taylor coefficients
of $a$ and $b$ are different.  Set
  $v:=(D\psi_x(x))^{-1}u$.
By \eqref{eq:expansion}, there exist
 an infinite subset
$I\subset\mathbb N$ and
a point $y\in\P^k$ such that, as $n\to\infty$ in $I$,
we have
\[
 x_n:=f_0^n(x)\longrightarrow y,
 \qquad \|D(f_0^n)_xv\|\longrightarrow\infty.
\]

Choose affine coordinates $\psi_y$ with $\psi_y(y)=0$.
Since  $K$ is compact,  there exists  $0<\rho<1$ such that
\[
 \sigma(\mathbb D_\rho)\subset\psi_y^{-1}\bigl(B(0,2)\bigr)
\]
whenever $\sigma\in K$ and
$\sigma(0)\in\psi_y^{-1}\bigl(\overline{B(0,1)}\bigr)$.
For all sufficiently large $n\in I$, both
$\F^n\gamma$ and $\F^n\eta$ are in this family, since their values at $0$ are $x_n$.

For each such $n$, the map
\[
 H_n(t,z):=\psi_y\bigl(f_t^n(\psi_x^{-1}(z))\bigr)
\]
is holomorphic near $(0,0)$. By Taylor expansion in $z$ at $b(t)$, we have 
\[
 H_n(t,b(t)+\zeta)
 =H_n(t,b(t))+D_zH_n(t,b(t))\zeta+O(\|\zeta\|^2).
\]
Since $b(t)=O(t)$ and $\delta(t)=t^r u+O(t^{r+1})$, we obtain
\begin{align*}
 H_n(t,a(t))-H_n(t,b(t))
 &=D_zH_n(t,b(t))\delta(t)+O(t^{2r})\\
 &=\bigl(D_zH_n(0,0)+O(t)\bigr)
       \bigl(t^r u+O(t^{r+1})\bigr)+O(t^{2r})\\
 &=t^rD_zH_n(0,0)u+O(t^{r+1}),
\end{align*}
where the last equality uses $2r\geq r+1$. By the chain rule, 
\begin{align*}
 [t^r]\bigl(\psi_y\circ\F^n\gamma-\psi_y\circ\F^n\eta\bigr)
 &=D_zH_n(0,0)u\\
 &=D\psi_y(x_n)\,D(f_0^n)_x(D\psi_x(x))^{-1}u\\
 &=D\psi_y(x_n)\,D(f_0^n)_xv,
\end{align*}
where $[t^r]$ denotes the coefficient of $t^r$.
Since $x_n\to y$, the Euclidean norm in the $\psi_y$ chart is
uniformly comparable to the Fubini--Study norm near $y$.
Thus the norm of this coefficient tends to infinity along $I$.
On the other hand, Cauchy's estimate on $|t|=\rho/2$ bounds it by
$4(2/\rho)^r$. This contradiction proves the proposition.
\end{proof}

\begin{cor}\label{cor:lifting}
	Let $\nu$ be an equilibrium web and let
	$K:=\supp\nu$. For every $\lambda_0\in M$, there is a Borel
	set $E\subset J_{\lambda_0}$ such that
	\[
	\mu_{\lambda_0}(E)=1,
	\qquad f_{\lambda_0}^{-1}(E)=E=f_{\lambda_0}(E),
	\]
	and evaluation induces a bijection
	\[
	e_{\lambda_0}:K_E\longrightarrow E,
	\qquad K_E:=\{\gamma\in K:\gamma(\lambda_0)\in E\}.
	\]
	Every point of $E$ has $d^k$ distinct preimages in $E$, and every
	inverse orbit in $E$ lifts uniquely to an inverse orbit in $K_E$.
\end{cor}

\begin{proof}
Set
	 $g=f_{\lambda_0}$ and $\mu=\mu_{\lambda_0}$. By \cite{BD99}, all Lyapunov exponents of $\mu$
	are positive. The multiplicative ergodic theorem therefore gives
	\eqref{eq:expansion} at $\mu$-almost every point.
	
	The equilibrium measure satisfies
	\[
	g_*\mu=\mu,\qquad g^*\mu=d^k\mu.
	\]
	It gives no mass to locally pluripolar sets, see
	\cite[Section~1]{DS10}. In particular, the critical set  has $\mu$-measure zero.  For every Borel set $A$, the preceding identities imply
	\begin{equation}\label{eq:null-sets}
		\mu(g^{-1}(A))=\mu(A),\qquad \mu(g(A))\leq d^k\mu(A).
	\end{equation}

	Choose a Borel null set $N\subset J_{\lambda_0}$ containing the
	critical  locus of $g$ in $J_{\lambda_0}$ and all points at which (\ref{eq:expansion}) fails.  Define the following Borel set
	\[
	E:=J_{\lambda_0}\setminus
	\bigcup_{a,b\geq0}g^{-a}\bigl(g^b(N)\bigr).
	\]
	
	Equation~\eqref{eq:null-sets} shows that $\mu(E)=1$.
	Complete invariance of $J_{\lambda_0}$ implies
	$g^{-1}(E)=E=g(E)$. All preimages of a point of $E$ lie in $E$
	and there are exactly $d^k$ of them.
	
	We have 
	\begin{equation}\label{eq:support}
		\F(K)=K,
		\qquad e_\lambda(K)=J_\lambda\quad(\lambda\in M).
	\end{equation}
 Thus $e_{\lambda_0}:K_E\to E$ is onto. By the definition of $E$,
condition~\eqref{eq:expansion} holds at every point $x\in E$.
Proposition~\ref{prop:uniqueness}, applied with this $x$, therefore
shows that at most one graph in $K$ takes the value $x$ at
$\lambda_0$. Hence $e_{\lambda_0}:K_E\to E$ is injective and thus
bijective.

Finally, let $(x_{-j})_{j\geq0}$ be an inverse orbit in $E$, and let
$\gamma_{-j}\in K_E$ be the unique graph satisfying
$\gamma_{-j}(\lambda_0)=x_{-j}$.
For every $j\geq0$, both $\F\gamma_{-j-1}$ and $\gamma_{-j}$
belong to $K$ and take the value $x_{-j}\in E$ at $\lambda_0$.
Since condition~\eqref{eq:expansion} holds at $x_{-j}$,
Proposition~\ref{prop:uniqueness}, applied with $x=x_{-j}$, gives
\[
\F\gamma_{-j-1}=\gamma_{-j}
\qquad\text{for every }j\geq0.
\]
Thus every inverse orbit in $E$ admits a unique lift to $K_E$.
\end{proof}

\section{Proof of Theorem \ref{thm:main}}\label{sec:stability}

	Let $\nu$ be an equilibrium web and set
	$K:=\supp\nu$. Fix $\lambda_0\in M$ and a repelling
	$n$-periodic point $p_0\in J_{\lambda_0}$. Replacing $f$ by $f^n$
	preserves the web and the equilibrium measures, so we may assume
	that $p_0$ is fixed.
 By \eqref{eq:support}, there is a holomorphic map
$p\in K$ such that
$p(\lambda_0)=p_0$.
Since $p_0$ is repelling, condition \eqref{eq:expansion} holds at
$p_0$. Proposition~\ref{prop:uniqueness} therefore shows that $p$ is
the unique element of $K$ taking the value $p_0$ at $\lambda_0$.
Since $\F p\in K$ and
\[
(\F p)(\lambda_0)=f_{\lambda_0}(p_0)=p_0,
\]
we obtain $\F p=p$.

	Choose a neighborhood $U$ of $p_0$ on which the inverse branch
	$\varphi$ of $f_{\lambda_0}$ fixing $p_0$ is contracting with $\varphi(U)\Subset U$. Let $E$ be the set
	given by Corollary~\ref{cor:lifting}. Set 
	\[
	A:=\{\gamma\in K_E:\gamma(\lambda_0)\in U\}.
	\]
	Since $p_0\in J_{\lambda_0}$ and $E$ has full measure,
	\begin{equation}\label{eq:positive-mass}
		\nu(A)=\mu_{\lambda_0}(U)>0.
	\end{equation}
	
	For each $\gamma\in A$, complete invariance of $E$ gives the inverse
	orbit
	\[
	z_{-j}:=\varphi^j(\gamma(\lambda_0))\in E,
	\qquad z_{-j}\longrightarrow p_0.
	\]
	By Corollary~\ref{cor:lifting},  there are unique holomorphic maps
	$\gamma_{-j}\in K_E$ satisfying
	\[
	\gamma_0=\gamma,\qquad
	\F\gamma_{-j-1}=\gamma_{-j},\qquad
	\gamma_{-j}(\lambda_0)=z_{-j}.
	\]
	Every subsequential limit of $(\gamma_{-j})$ in the compact set $K$
	takes the value $p_0$ at $\lambda_0$, so by 
	Proposition~\ref{prop:uniqueness}, the limit must be $p$.
	Consequently,
	\begin{equation}\label{eq:convergence}
		\gamma_{-j}\longrightarrow p
		\quad\text{locally uniformly on }M.
	\end{equation}
	
	Suppose that $p(\lambda)$ is not repelling for some $\lambda\in M$.
	The multiplier argument in \cite[Section~3]{BB26} gives a parameter
	$\lambda_1$ at which $p(\lambda_1)$ is hyperbolic and has at least
	one multiplier of modulus less than one. We may also arrange that
	$D f_{\lambda_1}(p(\lambda_1))$ is invertible.
	
Set
	  $g=f_{\lambda_1}$ and $q=p(\lambda_1)$. Choose a sufficiently
	small neighborhood $V$ of $q$ on which $g$ is injective, and let
	$W\subset V$ be its local unstable manifold at $q$. At least one
	multiplier is attracting, so $\dim_\C W<k$.
	For every $\gamma\in A$,  by equation~\eqref{eq:convergence} we have 
	\[
	\gamma_{-j}(\lambda_1)\longrightarrow q,
	\qquad
	g(\gamma_{-j-1}(\lambda_1))=\gamma_{-j}(\lambda_1).
	\]
	By the local
	unstable-manifold theorem, 
	\begin{equation}\label{eq:unstable-inclusion}
		\gamma(\lambda_1)\in\bigcup_{j\geq0}g^j(W)
		\quad\text{for every }\gamma\in A.
	\end{equation}

Since each $g^j(W)$ is locally pluripolar and
$\mu_{\lambda_1}$ gives no mass to locally pluripolar sets,
   we have 
	\[
	0<\nu(A)
	\leq\mu_{\lambda_1}\left(\bigcup_{j\geq0}g^j(W)\right)=0,
	\]
	a contradiction. Thus $p(\lambda)$ is $J$-repelling for every
	$\lambda\in M$, since $p\in K\subset\J$.

Applying the preceding construction to every point of
$R_n(\lambda_0)$, for every $n\geq1$, yields a finite family
$\mathcal R_n\subset\J$ of $n$-periodic graphs which are repelling
at every parameter. Distinct graphs in $\mathcal R_n$ are disjoint:
indeed, if two of them intersect at some parameter, then
Proposition~\ref{prop:uniqueness}, applied at that parameter, shows
that they coincide. Thus
\[
e_\lambda:\mathcal R_n\longrightarrow R_n(\lambda)
\]
is injective for every $\lambda\in M$.

Conversely, fix $\lambda\in M$ and $q\in R_n(\lambda)$. Repeating the
construction with $\lambda$ as the base parameter gives an
$n$-periodic graph in $K$ through $q$ which is repelling at every
parameter. Its value at $\lambda_0$ belongs to $R_n(\lambda_0)$, so
the graph belongs to $\mathcal R_n$. Hence $e_\lambda$ is also
surjective, and the proof is complete.

\begin{rem}\label{rem:compare-BBD}
For the applications in the preceding proof involving the repelling
periodic point $p_0$, the argument in the proof of
\cite[Lemma~2.5(3)]{BBD18} would be sufficient. In particular, that
argument gives the uniqueness of the periodic graph through $p_0$,
identifies the subsequential limits of the lifted graphs with this
periodic graph, and prevents two moving repelling periodic graphs
from intersecting.

The stronger Proposition~\ref{prop:uniqueness} is needed instead in
Corollary~\ref{cor:lifting}. The points of the completely invariant
set $E$ have full equilibrium measure but are not in general
periodic, so \cite[Lemma~2.5(3)]{BBD18} does not apply to them.
Although \cite[Lemma~2.5(2)]{BBD18} gives the existence of a graph in
$K$ through each such point, it gives neither uniqueness nor the
compatibility of these graphs under $\F$. Proposition
\ref{prop:uniqueness} provides precisely these two properties and
thereby allows a prescribed inverse orbit in $E$ to be lifted
coherently to an inverse orbit in $K$. This coherent lifting is what
is used above to obtain the convergence
\eqref{eq:convergence} for a positive-measure family of graphs.
\end{rem}

\section{Proof of Theorem \ref{thm}}

\begin{proof}
We first prove \textup{(i)}. Theorem~\ref{thm:main} gives
\textup{(C)}$\Rightarrow$\textup{(A)}, while
\cite[Corollary~1.7]{BBD18} gives
\textup{(A)}$\Rightarrow$\textup{(C')}. The implication
\textup{(C')}$\Rightarrow$\textup{(C)} is immediate. Thus
\textup{(A)}, \textup{(C)} and \textup{(C')} are equivalent.

The implication \textup{(A)}$\Rightarrow$\textup{(A')} is immediate.
Conversely, assume \textup{(A')}. After passing to a subsequence, we
may find integers $n_j\to\infty$ and $\tau>0$ such that
\[
d^{-kn_j}\#\mathcal P_{n_j}\longrightarrow\tau.
\]
Consider the finite measures
\[
\nu_j:=d^{-kn_j}
\sum_{\gamma\in\mathcal P_{n_j}}\delta_\gamma.
\]
Using the normality argument of
\cite[Proposition~2.2(2)]{BBD18}, we may extract a weak limit
$\nu$. The argument in the proof of
\cite[Theorem~4.1]{BB26} shows that
\[
(e_\lambda)_*\nu=\tau\mu_\lambda
\qquad\text{for every }\lambda\in M;
\]
here one uses the ergodicity of $\mu_\lambda$. Since the sets
$\mathcal P_{n_j}$ are $\mathcal F$-invariant, $\nu$ is also
$\mathcal F$-invariant. Hence $\tau^{-1}\nu$ is an equilibrium web,
and \textup{(A')}$\Rightarrow$\textup{(C)}.

Assume now \textup{(C')}. By
\cite[Proposition~2.4]{BBD18}, the family admits an acritical
ergodic equilibrium web. The proof of
\cite[Theorem~4.1]{BBD18}, which is local on the parameter space
once such a global web is given, produces an equilibrium lamination.
Thus \textup{(C')}$\Rightarrow$\textup{(D)} without any simple
connectedness assumption. Conversely,
\textup{(D)}$\Rightarrow$\textup{(A)} is
\cite[Proposition~1.2]{BB26}.
This proves \textup{(i)}.

We next prove \textup{(ii)}. The equivalence
\textup{(B)}$\Leftrightarrow$\textup{(E)} follows from
\cite[Proposition~3.12]{BBD18}. The same proposition, and more
precisely the mass estimate in
\cite[Lemma~3.13]{BBD18}, gives
\textup{(B)}$\Leftrightarrow$\textup{(E')}.
Moreover,
\textup{(B)}$\Leftrightarrow$\textup{(G)} follows from
\cite[Theorem~1.6]{BBD18}.
For the implication
\textup{(B)}$\Rightarrow$\textup{(G)}, see also
\cite[Theorem~A]{Bia19}, where it is proved more generally for
holomorphic families of polynomial-like maps of large topological
degree.
Finally,
\textup{(B)}$\Leftrightarrow$\textup{(F)} follows from
\cite[Theorem~2.2]{BBr23}, applied locally in the parameter space.
Moreover, \cite[Theorem~1.4]{BBr23} identifies
$\Omega(R_f)$ with the inverse image under $\pi_M$ of the local
$\mu$-stability locus. It follows that
\textup{(F)}$\Leftrightarrow$\textup{(F')}.
None of these implications requires simple connectedness.

For \textup{(iii)}, it is enough, in view of \textup{(i)}, to use
\textup{(A)}$\Rightarrow$\textup{(B)}, which is
\cite[Proposition~3.5]{BBD18}.

We now prove \textup{(iv)}. Assume that $M$ is simply connected.
Then
\textup{(B)}$\Rightarrow$\textup{(D)} follows from
\cite[Theorem~4.1]{BB26}. Together with \textup{(i)}--\textup{(iii)},
this proves that all the conditions are equivalent.

 Finally, to prove \textup{(v)}, assume that the equivalent conditions
in \textup{(i)} hold.
 By \textup{(A)}$\Rightarrow$\textup{(C')} and
\cite[Proposition~2.4]{BBD18}, the family admits an acritical
ergodic equilibrium web $\nu^{\rm eq}$. The uniqueness argument in
the proof of \cite[Theorem~4.1]{BBD18} is local on the parameter
space and compares $\nu^{\rm eq}$ with an arbitrary equilibrium web.
It therefore applies here without any simple connectedness
assumption and shows that $\nu^{\rm eq}$ is the unique equilibrium
web. The same argument gives
$\nu^{\rm eq}(\Lcal_1\mathbin{\triangle}\Lcal_2)=0$
for any two equilibrium laminations $\Lcal_1$ and $\Lcal_2$ of $f$.
\end{proof}

\begin{rem}\label{rem:new-uniqueness}
We observe that Proposition~\ref{prop:uniqueness} gives a direct alternative proof
of the uniqueness assertions in Theorem~\ref{thm}\textup{(v)}, using
only the existence of an equilibrium web.

Indeed, let $\nu_1$ and $\nu_2$ be two equilibrium webs and set
$K:=\operatorname{supp}\nu_1\cup\operatorname{supp}\nu_2$.
Then $K$ is compact and $\mathcal F(K)\subset K$. Fix
$\lambda_0\in M$. By the positivity of the Lyapunov exponents
\cite{BD99}, there is a Borel set $E\subset J_{\lambda_0}$ of full
$\mu_{\lambda_0}$-measure such that condition
\eqref{eq:expansion} holds at every $x\in E$. Hence
Proposition~\ref{prop:uniqueness} shows that the evaluation map
$e_{\lambda_0}$ is injective on
$K_E:=K\cap e_{\lambda_0}^{-1}(E)$.
Since
\[
(e_{\lambda_0})_*\nu_1=(e_{\lambda_0})_*\nu_2
=\mu_{\lambda_0}
\qquad\text{and}\qquad
\nu_1(K_E)=\nu_2(K_E)=1,
\]
it follows that $\nu_1=\nu_2$.

Similarly, if $\Lcal$ is an equilibrium lamination and
$\nu^{\rm eq}$ is the equilibrium web, we apply the proposition to
the compact forward-invariant set
$K:=\operatorname{supp}\nu^{\rm eq}\cup\overline{\Lcal}$.
Since $\mu_{\lambda_0}(e_{\lambda_0}(\Lcal))=1$, the same argument
gives $\nu^{\rm eq}(\Lcal)=1$. Applying this to two equilibrium
laminations $\Lcal_1$ and $\Lcal_2$ yields
$\nu^{\rm eq}(\Lcal_1\mathbin{\triangle}\Lcal_2)=0$.
\end{rem}

\begin{rem}
\label{rem:polynomial-like}
We briefly explain here
the polynomial-like extension stated after
Theorem~\ref{thm}. Let $d_t$ denote the topological degree. The
assumption of large topological degree means that $d_t$ is strictly
larger than the next dynamical degree $d^*_{k-1}$, and hence than all
the other dynamical degrees \cite{BDR24}.

The notation and the conditions in Theorem~\ref{thm} require the
following standard adaptations. Condition \textup{(A')} is understood
in the locally uniform form of \cite[Definition~1.3]{Bia19};
$d^{kn}$ is replaced by $d_t^n$; and the estimate in condition
\textup{(E')} is replaced, on every relatively compact parameter
set, by
$O(\delta^n)$
for some $\delta$ satisfying
$d_{k-1}^*<\delta<d_t$.
In conditions \textup{(F)} and \textup{(F')}, one uses the analogous
parameter--phase-space ramification current for holomorphic families
of polynomial-like maps, obtained from the Dinh--Sibony construction
\cite{DS03,DS10}, with the normalization $d_t^{-n}$ in place of
$d^{-kn}$. With this interpretation, the arguments proving
\textup{(B)}$\Leftrightarrow$\textup{(F)}
$\Leftrightarrow$\textup{(F')} in
\cite[Theorems~1.4 and~2.2]{BBr23} adapt to this setting.

The new implication
\textup{(C)}$\Rightarrow$\textup{(A)}
follows from the same proof as Theorem~\ref{thm:main}. The
polynomial-like counterparts of the other implications follow from
\cite{DS03,Pha05}, \cite[Chapters~2 and~3]{Bia16}, and
\cite{Bia19}. In particular,
\textup{(B)}$\Rightarrow$\textup{(G)}
is \cite[Theorem~A]{Bia19}, while the web--lamination statements and
the uniqueness of the equilibrium web are proved, in their original
simply connected formulation, in
\cite[Theorem~4.11]{Bia19}.
\end{rem}

\end{document}